\documentclass[12pt]{amsart}
\usepackage{amsmath,amsthm,amsfonts,amssymb, verbatim}
\usepackage{geometry}                
\usepackage{amssymb}
\usepackage{epstopdf}

\newcommand{\R}{\mathbb{R}}
\newcommand{\Z}{\mathbb{Z}}

\newcommand{\Pm}{\mathbb{P}}
\newcommand{\expE}{\mathbb{E}}
\newcommand{\no}{\nonumber}
\newcommand{\be}{\begin{equation}}
\newcommand{\ee}{\end{equation}}
\newcommand{\br}{\begin{eqnarray}}
\newcommand{\er}{\end{eqnarray}}
\newcommand{\bi}{\begin{itemize}}
\newcommand{\ei}{\end{itemize}}
\newcommand{\commentout}[1]{}
\newcommand{\norm}[1]{\lVert #1 \rVert}
\newcommand{\abs}[1]{\lvert #1 \rvert}
\newcommand{\eps}{\varepsilon}

\newcommand{\F}{\Psi}

\newtheorem{theorem}{Theorem}[section]
\newtheorem{lemma}[theorem]{Lemma}
\newtheorem{proposition}[theorem]{Proposition}
\newtheorem{corollary}[theorem]{Corollary}

\newtheorem*{theorem*}{Theorem}
\newtheorem*{lemma*}{Lemma}
\newtheorem*{proposition*}{Proposition}
\newtheorem*{corollary*}{Corollary}

\theoremstyle{definition}

\theoremstyle{remark}

\newtheorem*{remark*}{Remark}

\title[Homogenization of the G-equation with incompressible random drift]{Homogenization of the G-equation with incompressible random drift in two dimensions}
\author[James Nolen and Alexei Novikov]{James Nolen and Alexei Novikov}

\address{Department of Mathematics, Duke University, Durham, NC 27708, USA}
\email{nolen@math.duke.edu}
\address{Department of Mathematics, Penn State University,
University Park, State College, PA 16802}
\email{anovikov@math.psu.edu}

\thanks{The authors were partially supported by NSF grants DMS-1007572 and DMS-0908011, respectively.}

\dedicatory{To appear in Communications in Mathematical Sciences}

\begin{document}
\maketitle

\begin{abstract}
We study the homogenization limit of solutions to the G-equation with random drift. This Hamilton-Jacobi equation is a model for flame propagation in a turbulent fluid in the regime of thin flames. For a fluid velocity field that is statistically stationary and ergodic, we prove sufficient conditions for homogenization to hold with probability one. These conditions are expressed in terms of travel times for the associated control problem.  When the spatial dimension is equal to two and the fluid velocity is divergence-free, we verify that these conditions hold under suitable assumptions about the growth of the random stream function.
\end{abstract}

\section{Introduction}

We study the asymptotic behavior as $\epsilon \to 0$ of the solution of the initial value problem
\begin{gather}
u^\epsilon_t + V \left(\frac{x}{\epsilon},\omega \right) \cdot Du^\epsilon = \abs{Du^\epsilon} , \quad t > 0 , \;\; x \in \R^d, \label{Geqn} \\
u^\epsilon = u_0(x),\quad t = 0, \;\;x \in \R^d. \no
\end{gather}
The vector field $V:\R^d \times \Omega \to \R^d$ is assumed to be random, statistically stationary and ergodic with respect to $x$; $(\Omega, \mathcal{F}, \mathbb{P})$ is a given probability space, $\omega \in \Omega$. The initial condition $u_0(x)$ is assumed to be bounded and uniformly continuous. Under suitable hypotheses on $V$, we prove that with probability one the function $u^{\epsilon}$ converges (as $\epsilon \to 0$) locally uniformly in $[0,\infty) \times \R^d$ to a function $\bar u(t,x)$ which satisfies an equation of the form
\[
\bar u_t  = \bar H(D\bar u), \quad t > 0, \;\; x \in \R^d,
\]
with the same initial condition $\bar u(0,x) = u_0(x)$.  The function $\bar H:\R^d \to [0,\infty)$ is convex and homogeneous of degree one.

The level-set equation (\ref{Geqn}) is called the G-equation, and it is used as a model for turbulent combustion in the regime of thin flames \cite{Wlms85, Pet}. In this model, the level sets of $u^\epsilon$ represent the flame surface, and $V$ is the velocity of the underlying fluid (assumed to be independent of $u^\epsilon$). Wherever $u^\epsilon$ is differentiable and $\abs{Du^\epsilon} \neq 0$, the level sets of $u^\epsilon$ move with normal velocity $\nu = 1 + V \cdot \hat n$, where $\hat n = - Du^\epsilon/\abs{Du^\epsilon}$. When $V \equiv 0$, the level sets move with constant speed $s_L = 1$, which is called the laminar speed of flame propagation. An important scientific problem is to understand how turbulent transport affects the speed of the spreading flame \cite{Ro}. 

Let us state our assumptions on $V$. We suppose that $V:\R^d \times \Omega \to \R^d $ is a stationary random vector field, with $(\Omega, \mathcal{F}, \mathbb{P})$ a given probability space. Specifically, we assume that there is a group $\{\pi_x\}_{x \in \R^d}$ of transformations of $(\Omega, \mathcal{F})$ which is measure-preserving
\[
\Pm(\pi_x(A)) = \Pm(A), \quad \forall A \in \mathcal{F},
\]
and acts ergodically on $(\Omega,\mathcal{F})$ which means that if $A \in \mathcal{F}$ and $\pi_x(A) = A$ for all $x \in \R^d$, then either $\Pm(A) = 0$ or $\Pm(A) = 1$. We assume that the map $(x,\omega) \mapsto \pi_x \omega$ from $\R^d \times \Omega$ to $\R^d$ is jointly measurable, and that
\[
V(x,\omega) = \hat V(\pi_x \omega) , \;\; \forall x \in \R^d
\]
for some random variable $\hat V \in L^\infty(\Omega, \mathcal{F}, \mathbb{P})$.  We assume that with probability one, $V(\cdot,\omega)$ defines a vector field that is $C^1$, and satisfies the bound
\[
\norm{V(\cdot,\omega)}_{L^\infty(\R^d)} = \norm{\hat V}_{L^\infty(\Omega)} = V_\infty < \infty.
\]
We do not assume that $V_\infty < s_L = 1$, and this presents the fundamental mathematical difficulty, as we will describe. In particular, the Hamiltonian $H(p,x) = \abs{p} - V(x) \cdot p$ is not coercive in $p$ at every point $x \in \R^d$, so existing results on homogenization of random Hamilton-Jacobi equations (\cite{Soug1, RT}) do not apply here.

If the deterministic function $u_0$ is bounded and uniformly continuous, then with probability one, there exists a unique family of viscosity solutions $\{ u^\epsilon\}_{\epsilon > 0}$ to (\ref{Geqn}).  Our main result is a sufficient condition on $V$ which guarantees that homogenization occurs in dimension $d = 2$. In addition to the preceding assumptions on $V$, we assume that with probability one $V$ is divergence free: $\nabla \cdot V(x,\omega) = 0$ for all $x \in \R^2$.  In this case, the vector field $V$ may be associated with a stream function:  $V = \nabla^\perp \F = (- \partial_{x_2} \F, \partial_{x_1} \F)$, where $\F(x,\omega): \R^2 \times \Omega \to \R$.

\begin{theorem}\label{thm:twodim} Let $V:\R^2 \times \Omega \to \R^2$ be a stationary random vector field, uniformly bounded, $C^1$, and divergence-free, $\Pm$-almost surely. Suppose that $\expE[V(0,\omega)] = 0$ and that $V = \nabla^\perp \F = (- \partial_{x_2} \F, \partial_{x_1} \F)$ for some stream function $\F$ which satisfies
\be
\int_0^\infty \Pm \left( \sup_{\abs{x} \leq r} |\F(x,\omega) - \F(0,\omega) | > r/6 \right) \,dr < \infty. \label{streamgrowthcond}
\ee
Let $u_0(x)$ be bounded and uniformly continuous in $x \in \R^2$. Then there is a convex  function $\bar H:\R^2 \to [0,\infty)$ such that the following holds with probability one: for any $T > 0$ and $R > 0$,
\be
\lim_{\epsilon \to 0} \sup_{t \in [0,T]} \sup_{\abs{x} \leq R}  \left | u^\epsilon(t,x,\omega) - \bar u(t,x) \right| = 0 \label{homconv1}
\ee
where $\bar u$ is the unique viscosity solution of the initial value problem
\begin{gather}
\bar u_t = \bar H(D \bar u) \quad x \in \R^2,\; t > 0,  \\
\bar u(0,x) = u_0(x) , \quad x \in \R^2. \no
\end{gather}
The function $\bar H$ satisfies $\bar H(\lambda p) = \lambda \bar H(p)$ for all $\lambda > 0$ and $p \in \R^2$. Moreover, $\bar H(p) = 0$ if and only if $p = 0$. If the law of $V(x,\omega)$ is equal to the law of $V(\Phi x,\omega)$ for any orthogonal transformation $\Phi$, then $\bar H(p) \geq \abs{p}$ for all $p \in \R^2$.
\end{theorem}

\vspace{0.2in}
The starting point of our analysis is the control representation for the solution $u^\epsilon$:
\begin{equation}
u^\epsilon(x,t,\omega) = \sup_{\alpha \in \mathcal{A}_t} u_0(X^{\epsilon,\alpha}_x(t)). \label{controlrep}
\end{equation}
Here the set $\mathcal{A}_t$ is the set of controls 
\[
\mathcal{A}_t = \{ \alpha \in L^\infty([0,t];\R^d) \;:\; \| \alpha \|_\infty \leq 1 \},
\]
and for each $x \in \R^d$ and control $\alpha \in \mathcal{A}_t$, the function $X^{\epsilon,\alpha}_x(s)$ is defined by
\be
\frac{d}{ds} X^{\epsilon,\alpha}_x(s) = - V(\epsilon^{-1} X^{\epsilon,\alpha}_x(s)) + \alpha(s), \quad s \in [0,t], \quad \quad X^{\epsilon,\alpha}_x(0) = x. \label{controlode}
\ee
To keep the notation simple here, we have suppressed the dependence of $V$ and $X$ on $\omega \in \Omega$. See \cite{Barles} for proof of this representation.

To understand the behavior of $u^\epsilon$ as $\epsilon \to 0$, we will consider the random variables
\be
\tau(x,y,\omega) = \inf \{ t \geq 0\;|\; X^\alpha_x(t) = y, \;\;\;\text{for some}\; \alpha \in \mathcal{A}_t \}. \label{taudef1}
\ee
Here $\epsilon = 1$. We refer to $\tau(x,y,\omega)$ as the travel time from $x\in \R^d$ to $y\in \R^d$, since it is the first time that a path starting from $x$ may be controlled to the point $y$. From the stationarity of $V$ it follows that $\tau(x,y,\omega) = \tau(0,y-x,\pi_x \omega)$. Ultimately we will apply the subadditive ergodic theorem to this family of random variables to show that $r^{-1} \tau(0,rz,\omega)$ approaches a deterministic limit as $r \to \infty$. Then, through the control representation (\ref{controlrep}) this leads to the convergence of $u^\epsilon$ described in (\ref{homconv1}). However, because we allow $\abs{V}$ to be larger than $s_L = 1$, it is not clear a priori whether the travel time $\tau(x,y,\omega)$ is finite for every pair $(x,y)$.  Indeed, it is easy to construct vector fields for which  $\tau$ cannot be finite. For example, if $V =  \nabla Q(x)$ for some function $Q(x)$ satisfying $Q(x) = \abs{x}^2$ for $\abs{x} < 1$, then $\tau(0,y)$ cannot be finite for any $y$ with $\abs{y} > 1/2$. Because of the constraint $\abs{\alpha} \leq 1$, the control cannot overcome the strong flow $-V$ directed toward $x = 0$.  On the other hand, in Theorem \ref{thm:twodim} we are assuming that $V$ has zero divergence, so it is hopeful that this trapping phenomenon does not occur with such a velocity field; from each point $x$, there should always be a way of escape. 

In the case that $V$ varies periodically with respect to $x$ (and possibly time $t$), homogenization results analogous to (\ref{homconv1}) were proved recently by Cardaliaguet, Nolen, and Souganidis \cite{CNS} and by Xin and Yu \cite{XY}. As shown in \cite{CNS}, a sufficient condition for homogenization in all dimensions $d \geq 2$ is for $\norm{\nabla \cdot V}_{L^d(Q_1)}$ to be sufficiently small, where $Q_1$ is the period cell.  Homogenization may not hold, even in the periodic setting, for vector fields which do not satisfy the bound on $\nabla \cdot V$. An important part of these analyses is demonstrating the existence of ``approximate correctors"-- a family of periodic functions $\eta_{\lambda,p}(x)$ (with parameters $p \in \R^d$, $\lambda > 0$) which satisfy
\[
\lambda \eta_{\lambda,p}(x)  = \abs{p + D\eta_{\lambda,p}(x)} - V(x) \cdot (p + D\eta_{\lambda,p}(x)), \quad x \in \R^d.
\]
If one can obtain suitable bounds on the oscillation of $\lambda \eta_{\lambda,p}$, then $\lambda \eta_{\lambda, p} \to \bar H(p)$ as $\lambda \to 0$, and homogenization follows by standard arguments which extend the perturbed test-function method developed in \cite{Evans92}. For the G-equation, however, it is difficult to obtain such control because the Hamiltonian $H(p,x) = \abs{p} + V(x) \cdot p$ is not everywhere coercive in the gradient variable $p$ when we allow $\abs{V} > 1$. Homogenization with noncoercive Hamiltonians has been studied in other contexts \cite{AI01, ABMemoirs,  GBa, Ca, CLP, IM}, yet those results do not apply here, even in the periodic setting.

Beside the difficulties that arise from the noncoercivity of the Hamiltonian, the randomness of the medium presents challenges that are not present in periodic or almost periodic media. In particular, constructing and analyzing correctors or approximate correctors may be impossible. We do not pursue a precise discussion of this point here; the reference \cite{LSc} explores this issue in more detail. Also, in the random setting, the extreme behavior of the medium can lead to nonhomogenization \cite{EWX}. Instead of using correctors or approximate correctors, our approach to proving Theorem \ref{theo:homog} relies on the control representation (\ref{controlrep})-(\ref{controlode}). This is similar to the strategy used previously by Souganidis \cite{Soug1}, Rezakhanlou, and Tarver \cite{RT},  who proved homogenization of Hamilton-Jacobi equations in random media when the Hamiltonian is convex and coercive. In those works, it is the action functional that is analogous to our $\tau(x,y,\omega)$. Those results do not apply directly to the present setting, however, due to the noncoercivity of $H$. For the G equation in the periodic setting, Xin and Yu \cite{XY} used the control representation to obtain the necessary control of the approximate corrector $\eta_{\lambda,p}$.

The novelty of the present analysis is the method of controlling the travel times $\tau(x,y,\omega)$ so that the control representation (\ref{controlrep})-(\ref{controlode}) can be used to prove homogenization. As we have mentioned, we expect that some control on the divergence of $V$ is necessary. Even if we assume $\nabla \cdot V = 0$, however, other features of the medium may be obstacles to homogenization, and these obstacles may occur on arbitrarily large scales, since the medium is random. For example, there may be large vortices or shear layers which may require a large amount of time to traverse. In two dimensions, the geometry of $V$ can be controled through the stream function $\F$. Roughly speaking, the condition (\ref{streamgrowthcond}) precludes the frequent appearance of very large scale features that make $\tau(x,y,\omega)$ large. We do not know whether a condition like (\ref{streamgrowthcond}) is necessary for homogenization, or whether $\nabla \cdot V = 0$ alone is sufficient, as in the periodic case. If we make the additional assumption that the stream function $\F(x,\omega)$ is stationary, then condition (\ref{streamgrowthcond}) may be simplified:

\vspace{0.2in}

\begin{proposition} \label{prop:Fstat} Let $d = 2$. If $\F(x,\omega)$ is stationary, then condition (\ref{streamgrowthcond}) in Theorem \ref{thm:twodim} may be replaced by
\[
\expE \left[ \abs{\F(0)}^{3}  \right] < \infty. \label{Fmomentbound}
\]
\end{proposition}

\vspace{0.2in}

Our proof of Theorem \ref{thm:twodim} is based on the following theorem which gives sufficient conditions for homogenization in any dimension $d \geq 2$. For $R > 0$ we define the random variable
\[
\gamma(R,\omega) = \sup_{x,y \in B_R(0)} \tau(x,y,\omega).
\]

\begin{theorem}\label{theo:homog} Let $d \geq 2$, and let $u_0(x)$ be bounded and uniformly continuous in $x \in \R^d$. In addition to the aforementioned hypotheses on the random field $V$, assume that for all $R > 0$,
\be
\expE[\gamma(R,\omega)] < \infty, \label{gammaexp}
\ee
and that	
\be
\limsup_{R \to \infty}  \frac{1}{R} \gamma(R , \omega)    < \infty \label{taubound}
\ee
holds with probability one. Then there is a convex  function $\bar H:\R^d \to [0,\infty)$ such that the following holds with probability one: for any $T > 0$ and $R > 0$,
\be
\lim_{\epsilon \to 0} \sup_{t \in [0,T]} \sup_{\abs{x} \leq R}  \left | u^\epsilon(t,x,\omega) - \bar u(t,x) \right| = 0 \label{homconv}
\ee
where $\bar u$ is the unique viscosity solution of the initial value problem
\begin{gather}
\bar u_t = \bar H(D \bar u) \quad x \in \R^d,\; t > 0, \label{barueqn} \\
\bar u(0,x) = u_0(x) , \quad x \in \R^d. \no
\end{gather}
The function $\bar H$ satisfies $\bar H(\lambda p) = \lambda \bar H(p)$ for all $\lambda > 0$ and $p \in \R^d$. If $V$ is divergence-free, then $\bar H(p) = 0$ if and only if $p = 0$. If the law of $V(x,\omega)$ is equal to the law of $V(\Phi x,\omega)$ for any orthogonal transformation $\Phi$, then $\bar H(p) \geq \abs{p}$ for all $p \in \R^d$.
\end{theorem}

In the first part of Theorem \ref{theo:homog} we do not assume that $\nabla \cdot V = 0$. However, by the simple example described already, we know that conditions (\ref{gammaexp}) and (\ref{taubound}) may not hold without some restrictions on the divergence.  Since equation (\ref{Geqn}) is first order, information propagates at finite speed. We may think of $R/\gamma(R,\omega)$ as a lower bound on the speed at which information propagates between any two points in $B_R(0)$. Condition (\ref{gammaexp}) precludes this speed from being too small, on average. Condition (\ref{taubound}) controls extreme behavior of the random environment as one observes larger and larger regions of space. These two conditions may be difficult to verify, in general. However, in two dimensions and with the hypotheses of Theorem \ref{thm:twodim}, we can use the stream function $\F$ to prove that these conditions do hold.

Before proceeding to the proofs of these results, we mention that further analysis and numerical computation related to this equation with periodic drift may be found 
in~\cite{It, EMS, NXY, Ob_01} and references therein. In particular, the predictions of the G-equation model for flame propagation may differ significantly from the predictions of a reaction-diffusion-advection model \cite{EMS, NXY}. Some related results for reaction-diffusion-advection equations with random drift may be found in \cite{MS2, NXaihp}.

The rest of this article is organized as follows: In Section \ref{sec:ttasymp}, we will apply the subadditive ergodic theorem to show that the averaged travel times $\tau(x,y,\omega)/\abs{x-y}$ become deterministic as $\abs{x - y} \to \infty$. Eventually, we convert this into a statement about $u^\epsilon$, via the representation formula for $u^\epsilon$. However, that step requires some uniform control of $\tau(x,y,\omega)$, which we obtain in Lemma \ref{lem:tauuniform}. In Section \ref{sec:homog} we prove Theorem \ref{theo:homog} using the control representation for $u^\epsilon$ and the travel time estimates. The main step (Lemma \ref{lem:asympdense}) in this proof is to show that the domain of dependence for $u^\epsilon$ coincides approximately with that of $\bar u$. In Section \ref{sec:reachable} we prove a lower bound on the volume of the set of points that are reachable via the controlled trajectories (\ref{controlode}). That estimate plays an important role in Section \ref{sec:ttbound} where we prove Theorem \ref{thm:twodim} by verifying the conditions (\ref{controlrep})-(\ref{controlode}). Another key technical estimate in the proof of Theorem \ref{thm:twodim} is Lemma \ref{prop:travtimebound} which relates growth of the stream function $\F$ to an estimate on the travel times. Proposition \ref{prop:Fstat} is proved at the end of Section \ref{sec:ttbound}.

\section{Asymptotic behavior of the travel times} \label{sec:ttasymp}

In this section we describe the almost-sure asymptotic behavior of $\tau(x,y,\omega)$ as $\abs{x - y} \to \infty$. Except where stated otherwise, we assume thoughout this section that conditions (\ref{gammaexp}) and (\ref{taubound}) hold. The first result of this section is analogous to a ``shape theorem" for the first-passage time in percolation theory:

\begin{lemma} \label{lem:subadd}
There is a function $\bar q(p):\R^d \to [0,\infty)$ and a measurable set $\bar \Omega \subset \Omega$ such that $\Pm(\bar \Omega) = 1$ and 
\[
\lim_{r \to \infty} \frac{1}{r} \tau(0,r p,\omega) = \bar q(p)
\]
holds for all $p \in \R^d$ and $\omega \in \bar \Omega$. This function $\bar q$ is continuous, convex, positively homogeneous of degree one, and $\bar q(p) \geq \abs{p}/(1 + V_\infty)$ for all $p \in \R^d$. 
\end{lemma}

Due to the stationarity of $V$, this lemma implies that for any $x$ and $p$, $\tau(rx,r(x + p),\omega)/r \to \bar q(p)$ in probability as $r \to \infty$. However, to apply this result in the control representation for $u^\epsilon$, we will require some uniform control over $\tau(x,y,\omega)$, as described by the next Lemma:

\begin{lemma} \label{lem:tauuniform}  Let $\bar q$ and $\bar \Omega$ be as in Lemma \ref{lem:subadd}. There is a measurable set $\tilde \Omega \subset \bar \Omega$ with $\Pm(\tilde \Omega) = 1$ such that for each $R > 0$, $M > 0$, 
\be
\lim_{r \to \infty} \sup_{\abs{x} \leq R} \sup_{\substack{p \in \R^d\\ \abs{p} \leq M}}  \left| \frac{\tau(rx,r(x + p),\omega) }{r} - \bar q(p) \right| = 0. \label{Gdeltconv}
\ee
holds for every $\omega \in \tilde \Omega$.
\end{lemma}

In proving Lemma \ref{lem:subadd} and Lemma \ref{lem:tauuniform}, we will make use of the following consequence of condition (\ref{taubound}). We give a proof of this proposition at the end of this section. The proof does not require the condition (\ref{gammaexp}) to hold.

\begin{proposition}\label{lem:taubound2} Suppose that condition (\ref{taubound}) holds. Then 
\be
\lim_{\epsilon \to 0} \limsup_{r \to \infty}   \left( \sup_{ \abs{x} \leq r } \frac{1}{r} \gamma(r \epsilon, \pi_x \omega)  \right)  = 0\label{taubound2}
\ee
holds with probability one.
\end{proposition}

{\bf Proof of Lemma \ref{lem:subadd}:} First, we fix a vector $p \in \R^d$, $p \neq 0$ and define a non-negative family of random variables 
\[
q_{j,k}(p,\omega) = \tau(j p,k p,\omega), \quad 0 \leq j \leq k, \quad j,k \in \R.
\]
From the definition of $\tau$, it is easy to see that for any three points $x,y,z \in \R^d$, the triangle inequality holds for $\tau$:
\[
\tau(x,z,\omega) \leq \tau(x,y,\omega) + \tau(y,z,\omega).
\]
Consequently, $q_{j,k}$ is a subadditive family:
\[
q_{j,l} \leq q_{j,k} + q_{k,l} 
\]
for all $0 \leq j \leq k \leq l$. Moreover from the stationarity of $V$ it follows that
\[
q_{j+m,k+m}(\omega) = q_{j,k}(\pi_{m p} \omega)
\]
holds for all indices $j \leq k$, and all $m \geq 0$, so this is a stationary process. Finally, condition (\ref{gammaexp}) and the subadditivity imply that
\[
0 \leq \expE[q_{0,k} ]  \leq k C
\]
holds for all $k \geq 0$ for a constant $C$ that depends on $\abs{p}$. Therefore, by the subadditive ergodic theorem (for example, \cite{Lig85}) there is a random variable $\tilde q(p,\omega)$ such that the limit along integer values $k$
\be
\lim_{\substack{k \to \infty } } \frac{1}{k} q_{0,k}(p,\omega) = \tilde q(p,\omega) \label{qlim}
\ee
holds with with probability one.  

In fact, there is a constant $\bar q(p)$ such that $\tilde q(p,\omega) = \bar q(p)$ holds with probability one. This follows from the fact that the limit $\tilde q$ is invariant under $\pi_x$, as we now show. Suppose $x \in \R^d$ with $\abs{x} \leq R$.
\br
\tilde q(p,\pi_x \omega) & = & \lim_{k \to \infty} \frac{1}{k} \tau(0,kp,\pi_x \omega) \no \\
& = &  \lim_{k \to \infty} \frac{1}{k} \left( \tau(0,kp,\omega) + \tau(0,kp,\pi_x \omega)  -  \tau(0,kp,\omega) \right) \no \\
& = &  \lim_{k \to \infty} \frac{1}{k} \left( \tau(0,kp,\omega) + \tau(x,x + kp, \omega)  -  \tau(0,kp,\omega) \right) \label{qshift}
\er 
From the definition of $\gamma(R,\omega)$, we see that
\[
\abs{\tau(x,x + kp, \omega)  -  \tau(0,kp,\omega) } \leq \gamma(R,\omega) + \gamma(R,\pi_{kp} \omega).
\]
By the ergodic theorem, the limit
\[
\lim_{k\to \infty} \frac{1}{k} \sum_{n=1}^k \gamma(R,\pi_{np} \omega) 
\]
exists with probability one (the subset of $\Omega$ on which the limit exists may depend on $p$ and $R$, but not on $x$), and it is finite since $\expE[\gamma(R, \omega) ] < \infty$. Consequently,
\be
\lim_{k\to \infty} \frac{1}{k} \gamma(R,\pi_{kp} \omega) = \lim_{k\to \infty} \left( \frac{1}{k} \sum_{n=1}^k \gamma(R,\pi_{np} \omega) - \frac{1}{k} \sum_{n=1}^{k-1} \gamma(R,\pi_{np} \omega) \right)= 0 \label{gammalim}
\ee
holds with probability one.  Therefore, with probability one, 
\[
\lim_{k \to \infty} \frac{1}{k} \left( \tau(x,x + kp, \omega)  -  \tau(0,kp,\omega)\right) = 0
\]
holds for all $x \in \R^d$ with $\abs{x} \leq R$. Since $R$ is arbitrary, we conclude from (\ref{qshift}) that with probability one,
\[
\tilde q(p,\pi_x \omega) = \lim_{k \to \infty} \frac{1}{k} \tau(0,kp,\pi_x \omega) = \lim_{k \to \infty} \frac{1}{k} \tau(0,kp,\omega)  = \tilde q(p,\omega)
\]
holds for all $x \in \R^d$. Now the assumption that $\pi_x$ is an ergodic transformation implies that $\tilde q(p,\omega) = \bar q(p)$ with probability one.

We claim that limit (\ref{qlim}) holds along continuous time $k \in \R$. To see this, suppose that $r \in [n,n+1)$ where $n \geq 0$ is an integer. We have
\[
\tau(0,r p,\omega)  \leq  \tau(0,np,\omega) + \tau(np,rp,\omega) \leq \tau(0,np,\omega) + \gamma(\abs{p},\pi_{np} \omega). 
\]
Now using (\ref{gammalim}) with $R = \abs{p}$, we conclude that, with probability one,
\be
\lim_{n \to \infty} \sup_{r \in [n,n+1)} \frac{1}{r} \tau(0,r p,\omega) \leq \lim_{n \to \infty} \frac{1}{n} \tau(0,np,\omega) + \frac{1}{n}\gamma(\abs{p},\pi_{np} \omega) = \bar q(p) 
\ee
This, and a similar lower bound establishes the claim.

We have shown that for each $p \in \R^d$, there is a set $\Omega_p$ with $\Pm(\Omega_p) = 1$ on which (\ref{qlim}) holds along continuous time $k \in \R$. Now take $\bar \Omega = \cap_{p \in \mathbb{Q}^d} \Omega_p$ to obtain a measurable set $\bar \Omega \subset \Omega$ such that $\Pm(\bar \Omega = 1)$ and for all $\omega \in \bar \Omega$, the limit
\be
\lim_{r \to \infty } \frac{1}{r} q_{0,r}(\omega) = \bar q(p)  \label{qlimc}
\ee
holds along continuous time $r \in \R$, for all rational vectors $p \in \mathbb{Q}^d$.

From the definition of $\tau$ and the fact that $\abs{V}$ is bounded by $V_\infty$, it is clear that $\bar q(p) \geq \abs{p}/(1 + V_\infty)$ for all $p$. Using the subadditivity property, it is not hard to show that $\bar q(p)$ is continuous and convex in $p$ and homogeneous of degree one: $\bar q(\lambda p) = \lambda \bar q(p)$ for all $\lambda > 0$. Thus we may extend the definition of $\bar q(p)$ to all $p \in \R^d$. Let us now show that (\ref{qlimc}) holds for all $p \in \R^d$, not just for vectors $p \in \mathbb{Q}^d$. Suppose that $p_1 \in \R^d$. Let $p_2 \in \mathbb{Q}^d$ with $\abs{p_1 - p_2} \leq \epsilon$. Then
\br
\tau(0,r p_1,\omega)  \leq  \tau(0,r p_2, \omega) + \tau(r p_2, r p_1,\omega)  \leq \tau(0,r p_2, \omega) + \gamma(r \epsilon, \pi_{r p_2} \omega) \label{tbapp1}
\er 		
holds for all $\omega \in \bar \Omega$. Similarly, $\tau(0,r p_2,\omega)  \leq  \tau(0,r p_1, \omega) + \tau(r p_1, r p_2,\omega)$ so that
\be
\tau(0,r p_1, \omega)  \geq  \tau(0,r p_2,\omega) - \gamma(r \epsilon, \pi_{r p_2} \omega). \label{tbapp2}
\ee

Assumption (\ref{taubound}) and Proposition \ref{lem:taubound2} imply that for any $M > 0$
\be
\lim_{\epsilon \to 0} \limsup_{r \to \infty}   \left( \sup_{\abs{p} \leq M} \frac{1}{r} \gamma(r \epsilon, \pi_{r p} \omega) \right) = 0. 
\ee
holds with probability one.  By removing from $\bar \Omega$ a set of measure zero, if necessary, we may apply this to (\ref{tbapp1}) and (\ref{tbapp2}) and conclude that for some function $\nu(\epsilon) \geq 0$ which is $o(1)$ as $\epsilon \to 0$,
\be
\bar q(p_2)- \nu(\epsilon) \leq\liminf_{r \to \infty} \frac{1}{r} \tau(0,r p_1,\omega) \leq \limsup_{r \to \infty} \frac{1}{r} \tau(0,r p_1,\omega)  \leq \bar q(p_2) + \nu(\epsilon) \label{qp2upperlower}
\ee
holds for all $\omega \in \bar \Omega$. Since $\bar q$ is continuous and since $\epsilon$ may be made arbitrarily small by choosing $p_2 \in \mathbb{Q}^d$ arbitrarily close to $p_1$, this completes the proof of Lemma \ref{lem:subadd}.
\qed

In order to prove Lemma \ref{lem:tauuniform}, we will employ some ideas that have been used in \cite{KRV, KosVar, Sch}.

\begin{lemma}\label{lem:Neta} Let $\bar \Omega$ be as in Lemma \ref{lem:subadd}. Fix $\alpha > 0$. For each $\eta > 0$ there is a measurable set $N_\eta \subset \bar \Omega$ such that $\Pm(N_\eta) \geq 1 - \eta$ and 
\be
\lim_{r \to \infty} \sup_{\omega \in N_\eta} \sup_{\substack{p \in \R^d\\ \abs{p} \leq \alpha}}  \left| \frac{\tau(0,rp,\omega) }{r} - \bar q(p) \right| = 0. \label{Ndeltconv}
\ee
\end{lemma}
{\bf Proof of Lemma \ref{lem:Neta}:} Egorov's theorem implies that for each $p \in \R^d$ and $\eta > 0$, there is a measurable set $N_\eta^p \subset \bar \Omega$ such that $\Pm(N_\eta^p) \geq 1 - \eta$ and
\[
\lim_{r \to \infty} \sup_{\omega \in N_\eta^p} \left| \frac{\tau(0,rp,\omega) }{r} - \bar q(p) \right| = 0
\]
holds (i.e. uniformly over $N_\eta^p$).  For $\mathbb{Q}^d = \{ p_n \}_{n=1}^\infty$, let 
\[
N_\eta  = \bigcap_{n=1}^\infty N^{p_n}_{\eta 2^{-n}}
\]
Then $\Pm(N_\eta) \geq 1 - \eta$, and 
\be
\lim_{r \to \infty} \sup_{\omega \in N_\eta} \left| \frac{\tau(0,rp,\omega) }{r} - \bar q(p) \right| = 0 \label{Ndeltconv2}
\ee
holds for all rational vectors $p \in \mathbb{Q}^d$. The locally uniform convergence described by (\ref{Ndeltconv}) now follows from (\ref{Ndeltconv2}) and condition (\ref{taubound}) and Proposition \ref{lem:taubound2}, as in the derivation of (\ref{qp2upperlower}). \qed

\vspace{0.2in}

\begin{lemma} \label{lem:xxp} Let $\bar \Omega$ be as in Lemma \ref{lem:subadd}. There is a measurable set $\tilde \Omega \subset \bar \Omega$ with $\Pm(\tilde \Omega) = 1$ such that the following holds: For every $R > 0$, every $\omega \in \tilde \Omega$, and every integer $n$ sufficiently large (depending only on $d$), there exists a constant $r_0 = r_0(R,n, \omega)$ such that for all $r \geq r_0$ and $x \in \R^d$ satisfying $\abs{x} \leq R $, there exists at least one point $x' \in \R^d$ (depending on $r$ and $\omega$) satisfying
\[
\abs{x - x'} \leq  R \left(\frac{3}{2^n} \right)^{1/d},
\]
and
\[
\pi_{r x'} \omega \in N_{2^{-n}}.
\]
where $N_{2^{-n}}$ is defined in Lemma \ref{lem:Neta} with $\eta = 2^{-n}$.
\end{lemma}
{\bf Proof of Lemma \ref{lem:xxp}:} This may be proved as in \cite{KRV}; see the proof of Theorem 2.1 therein. Also, see Lemma 5.7 of \cite{Sch}. We provide here a  proof for reader's convenience.

Let $R>0$ and $\eta=2^{-n}$. By the ergodic theorem, there is a set of full measure $\Omega_R \subset \bar \Omega$ such that
\[
\lim_{r \to \infty} \frac{\abs{\{ y \in \R^d \;|\; \pi_{y} \omega \in N_\eta, \;\; \abs{y} \leq Rr \}}}{\abs{\{ y \in \R^d \;|\; \abs{y} \leq Rr \} } } = \Pm(N_\eta) \geq 1 - \eta
\]
holds for all $\omega \in \Omega_R$. Consequently,  there is $r_0 = r_0(R,n,\omega)$ such that for all $r > r_0$ and $h \geq 0$,
\br
\abs{ \{ z \in \R^d \;|\; \pi_{z} \omega \notin N_\eta, \;\; \abs{z} \leq R(r + h) \} } & \leq &  2 \eta \abs{\{ z \in \R^d \;|\; \abs{z} \leq R(r+h) \} } \no \\
& = & 2 \eta \left(1 + \frac{h}{r} \right)^d \abs{\{ z \in \R^d \;|\; \abs{z} \leq Rr \} }. \label{Rrupper}
\er

Let $h = (3 \eta)^{1/d} r$. For any $y \in B_{rR}(0)$, $B_\delta(y) \subset B_{(r + h)R}(0)$ if $\delta = hR$. Also
\be
\abs{B_{\delta}(y)} = \frac{\delta^d}{(Rr)^d} \abs{B_{rR}(0)}= \frac{\delta^d}{(Rr)^d} \abs{\{ z \;|\; \abs{z} \leq Rr \} }  = 3 \eta \abs{\{ z \;|\; \abs{z} \leq Rr \} }. \label{Rrlower}
\ee
If $\eta < ((5/4)^{1/d} - 1)^d/3$, we have for all $r> 0$, $2 \eta \left( 1 + \frac{h}{r}\right)^d	=  2 \eta \left( 1 + (3 \eta)^{1/d} \right)^d < 5 \eta /2$. Consequently, from (\ref{Rrupper}) and (\ref{Rrlower}) we see that if $r > r_0(R,n,\omega)$ there cannot be $y \in B_{rR}(0)$ such that 
\[
B_{\delta} (y) \subset \{ z \;|\; \pi_{z} \omega \notin N_\eta, \;\; \abs{z} \leq R(r + h) \}.
\]
So, for every $y \in B_{rR}(0)$ and $\omega \in \Omega_R$ there must be a point 
$y'$ such that $\abs{y - y'} \leq \delta$ and $\pi_{y'} \omega \in N_\eta$. 

We finally set
$\tilde \Omega = \cap_{k = 1}^\infty \Omega_k$,
and the result follows with $x = y/r$ and $x' = y'/r$.
\qed

\vspace{0.2in}

{\bf Proof of Lemma \ref{lem:tauuniform}:} Let $\tilde \Omega$ be as in Lemma \ref{lem:xxp}. Let $\omega \in \tilde \Omega$. For any pair of points $x, x' \in \R^d$.
\br
\tau(rx,r(x + p),\omega) & = & \tau(r x',r(x'  + p),\omega) + \tau(r  x,r( x + p),\omega) - \tau(r x',r(x' + p),\omega) \no \\
&  = & \tau(0,r p,\pi_{r x'} \omega) + \tau(r  x,r( x + p),\omega) - \tau(r x',r(x' + p),\omega)
\er
so that
\br
\left| \frac{\tau(rx,r(x + p),\omega) }{r} - \bar q(p) \right| & \leq & \left| \frac{\tau(0,r p,\pi_{r x'} \omega) }{r} - \bar q(p) \right| \no \\
& &  + \frac{1}{r} \abs{\tau(r  x,r( x + p),\omega) - \tau(r x',r(x' + p),\omega)} \no \\
& \leq & \left| \frac{\tau(0,r p,\pi_{r x'} \omega) }{r} - \bar q(p) \right| \no \\
& &  + \frac{1}{r} \gamma(r \abs{x' - x}, \pi_{rx} \omega) +  \frac{1}{r} \gamma(r \abs{x' - x}, \pi_{r(x + p)} \omega) \label{tterms}
\er
Fix $n \in \mathbb{N}$ and $\abs{x} \leq R$, and let $r_0 = r_0(R,n,\omega)$ be as in Lemma \ref{lem:xxp}. Then for $r > r_0$ and $\abs{x} \leq R$, we may choose $x' = x'(x,r,n,\omega)$ according to Lemma \ref{lem:xxp}. Since $\pi_{r x'} \omega \in N_{2^{-n}}$, we have
\be
\sup_{\substack{p \in \R^d\\ \abs{p} \leq \alpha}}  \left| \frac{\tau(0,r p,\pi_{r x'} \omega) }{r} - \bar q(p) \right| \leq \sup_{\omega' \in N_{2^{-n}}} \sup_{\substack{p \in \R^d\\ \abs{p} \leq \alpha}}  \left| \frac{\tau(0,r p,\omega') }{r} - \bar q(p) \right|  \label{xxpbound3}
\ee
and the right hand side is independent of $x$, for all $\abs{x} \leq R$. To control the other terms in (\ref{tterms}) we may apply (\ref{taubound}) and (\ref{taubound2}), since $\abs{x - x'} \leq R(\frac{3}{2^n})^{1/d}$ may be made arbitrarily small by taking $r$ and $n$ large, with $R$ fixed. Consequently, (\ref{taubound}) and (\ref{taubound2}) imply that
\begin{gather}
\lim_{n \to \infty}  \limsup_{r \to \infty} \sup_{\abs{x} \leq R} \frac{1}{r} \gamma(r \abs{x' - x}, \pi_{rx} \omega) = 0 \label{xxpbound1} 
\end{gather}
and
\begin{gather}
\lim_{n \to \infty}  \limsup_{r \to \infty} \sup_{\abs{x} \leq R} \sup_{\abs{p} \leq \alpha} \frac{1}{r} \gamma(r \abs{x' - x}, \pi_{r(x + p)} \omega) = 0. \label{xxpbound2} 
\end{gather}

Because $n \in \mathbb{Z}$ may be chosen arbitrarily large as $r \to \infty$, we combine (\ref{tterms}),  (\ref{xxpbound3}), (\ref{xxpbound1}), (\ref{xxpbound2}), and Lemma \ref{lem:Neta} to conclude that for all $\omega \in \tilde \Omega$, $R > 0$, and $\alpha > 0$,
\be
\limsup_{r \to \infty} \sup_{\abs{x} \leq R} \sup_{\substack{p \in \R^d\\ \abs{p} \leq \alpha}}  \left| \frac{\tau(rx,r(x + p),\omega) }{r} - \bar q(p) \right| = 0.
\ee
This proves Lemma \ref{lem:tauuniform}. \qed

\vspace{0.2in}
{\bf Proof of Proposition \ref{lem:taubound2}:} Assume that condition condition (\ref{taubound}) holds (but we need not assume condition (\ref{gammaexp})). For $h > 0$ and $M_0 > 0$, let $G_{M_0,h} \subset \Omega$ denote the set 
\[
G_{M_0,h} = \left \{ \omega \in \Omega \;|\;\; \frac{1}{M} \gamma(M,\omega) \leq h,\;\; \forall M \geq M_0 \right\}.
\]
By condition (\ref{taubound}) we know that for any $\epsilon > 0$ we may choose $h$ sufficiently large and $M_0$ sufficiently large so that $\Pm(G_{M_0,h}) \geq 1 - \epsilon$. Observe that for all $y \in \R^2$, $\Pm(\pi_y G_{M_0,h}) = \Pm(G_{M_0,h})$. Therefore, the ergodic theorem implies that
\be
\lim_{R \to \infty} \frac{\left| \{ x \in B_R(0) \;|\; \pi_x \omega \in G_{M_0,h} \} \right|}{\abs{B_R(0)}} = \Pm\left( G_{M_0,h}\right) 
\ee
holds with probability one. So, for almost every $\omega \in \Omega$, there is an $R_0 = R_0(M_0,h,\omega)$ such that
\[
\left| \{ x \in B_R(0) \;|\; \pi_x \omega \notin G_{M_0,h} \} \right| \leq 2 \abs{B_R(0)} \left(1 - \Pm( G_{M_0,h}) \right) 
\]
holds for all $R \geq R_0$. Thus, if $R \geq R_0$ and $x \in B_R(0)$, there must be a point $x' \in B_R(0)$ such that $\abs{x - x'} \leq C R (1 - \Pm( G_{M_0,h}))^{1/d}$ and $\pi_{x'} \omega \in G_{M_0,h}$, where $C$ is a universal constant. 

Now, for any $\epsilon > 0$, we may choose $h$ and $M_0$ sufficiently large so that $C  (1 - \Pm( G_{M_0,h}))^{1/d} \leq \epsilon$. Thus, for $R \geq R_0(M_0,h,\omega)$, every point $x \in B_R(0)$ is contained in a ball of the of form $B_{R \epsilon}(x')$ where $\pi_{x'} \omega \in G_{M_0,h}$. Since the point $x' \in G_{M_0,h}$ we know that
\be
\gamma(M,\pi_{x'} \omega) \leq h M, \quad \forall \;M \geq M_0.
\ee
In particular, $\gamma(2 R \epsilon,\pi_{x'} \omega) \leq 2 h R \epsilon$, for $R \geq M_0 \epsilon^{-1}/2$. On the other hand, $B_{2R\epsilon}(x') \supset B_{R\epsilon}(x)$ since $\abs{x - x'} \leq R\epsilon$, so the definition of $\gamma$ now implies that
\[
\gamma(R \epsilon,\pi_x \omega) \leq \gamma(2 R \epsilon,\pi_{x'} \omega) \leq 2 h R \epsilon
\]
holds for all $x \in B_R(0)$ and $R \geq \max(R_0, M_0 \epsilon^{-1}/2)$. Taking $R \to \infty$ we conclude that
\[
\limsup_{R \to \infty} \frac{1}{R} \sup_{\abs{x} \leq R} \gamma(R\epsilon,\pi_x \omega) \leq 2 h \epsilon
\]
holds with probability one. \qed

\section{Homogenization} \label{sec:homog}
In this section we prove Theorem \ref{theo:homog} using the estimates on the travel times (Lemma \ref{lem:tauuniform}) and the control representation for $u^\epsilon(t,x,\omega)$. Throughout this section we assume that conditions (\ref{gammaexp}) and (\ref{taubound}) hold. First, we identify the effective Hamiltonian for which Theorem \ref{theo:homog} holds. Let
\[
\bar H(p) = \sup \left \{ p \cdot z \;|\;\; z \in \R^d, \;\; \bar q(z) = 1 \right \}.
\]
Since $\bar H(p)$ is the supremum of a family of linear functions of $p$, it is immediate that $\bar H$ is convex in $p$, and positively homogeneous of degree one. Since $\bar q(p) \geq \abs{p}/(1 + V_\infty)$ it follows that $\bar H(p) \leq \abs{p}(1 + V^\infty)$.

Fix $\omega \in \tilde \Omega$ where $\tilde \Omega$  as in Lemma~\ref{lem:tauuniform}. For any $x \in \R^d$ and $t > 0$ we define $\Gamma_{x,t}(\omega) \subset \R^d$ to be the set of points $y$ for which $X^\alpha_x(t) = y$ for some control $\alpha \in \mathcal{A}_t$:
\be\label{reachable}
\Gamma_{x,t}(\omega)=\{ y \in \mathbb{R}^d\;|\;X^{\alpha}_{x}(t) = y, \quad \text{for some}\; \alpha \in {\mathcal{A}_t} \},
\ee	
and $X^{\alpha}_{x}(t)$ solves~\eqref{controlode} with $\eps=1$.
We refer to this set as the reachable set at time $t$, starting from $x$. With $\epsilon = 1$, the control representation for $u$ is
\be
u^1(x,t,\omega) = \sup_{\alpha \in \mathcal{A}_t} u_0(X^\alpha_x(t)) = \sup_{y \in \Gamma_{x,t}(\omega)} u_0(y). \label{contrp}
\ee
Thus, the set $\Gamma_{x,t}$ may be regarded as the domain of dependence of the solution $u^1(x,t)$. We will compare this set with the bounded, convex set $\{x + W_t \} = \{ x + y \;|\; y \in W_t \}$ where
\[
W_t = \{ t v \;|\;\; v \in \R^d,\;\; \bar q(v) \leq 1 \} = \{ z \;|\;\; z \in \R^d,\;\; \bar q(z) \leq t \}. 
\]
The set $\{x + W_t\}$ is the domain of dependence of the function $\bar u(x,t)$ satisfying the effective equation. The following lemma shows that for $t$ large, the reachable set $\Gamma_{x,t}$ approximately coincides with $\{x + W_t\}$:

\begin{lemma}\label{lem:asympdense} Let $\delta > 0$, $R > 0$. For all $\omega \in \tilde \Omega$, there is a time $t_0 = t_0(\delta,R,\omega) > 1$ such that if $t \geq t_0$, then the following hold:
\bi
\item[(i)] For all $\abs{x} \leq Rt$, $\Gamma_{x,t}(\omega) \subset \{ x + W_{t(1 + \delta)}\}$.
\item[(ii)] For all $\abs{x} \leq Rt$ and all $z \in \{ x + W_t\}$ there is a point $y \in \Gamma_{x,t}(\omega)$ such that $\abs{y - z} \leq \delta t$.
\ei
\end{lemma}
{\bf Proof of Lemma \ref{lem:asympdense}:} First we prove (i). Observe that
\[
\Gamma_{x,t}(\omega) \subset \{ \hat x \in \R^d \;|\; \tau(x,\hat x,\omega) \leq t \}.
\]
Let $m = 2(\norm{V}_\infty + 1)$. From the definition of $X^\alpha(t)$ it follows that if $\abs{v} \geq m$, then $\tau(x,x + tv) \geq 2t$ for all $x \in \R^d$ and $t > 0$. Therefore,
\[
\Gamma_{x,t}(\omega)  \subset \{ \hat x \;|\; \tau(x,\hat x,\omega) \leq t \} \subset  \{ x + tv  \;|\; \abs{v} \leq m, \;\; \tau(x,x + tv,\omega) \leq t \}
\]
holds for all $t > 0$. Now we apply Lemma \ref{lem:tauuniform} to conclude that for $t$ sufficiently large (depending on $\delta$, $R$, and $\omega$, but not on $x$)
\br
\Gamma_{x,t}(\omega) & \subset & \left \{ x + t v \in \R^d \; | \; \bar q(v) \leq 1 + \delta\right\} 
\er
holds for all $\abs{x} \leq Rt$. This last set is precisely $\left \{ x + t v \in \R^d \; | \; \bar q(v) \leq 1 + \delta\right\}  = \{ x + W_{t(1 + \delta)}\}$. This proves that there is $t_0 = t_0(\delta,R,\omega)$ such that (i) holds for all $t \geq t_0$.

Now we prove (ii). We first prove that (ii) holds for $z$ in the boundary of the set $\{ x + W_t\}$. Then we will prove it for $z$ in the interior of this set. Let 
\[
\mathcal{H} = \partial W_1 = \{ z \;|\;\; z \in \R^d,\;\; \bar q(z) = 1 \},
\]
and let $M > 1 + V_\infty$ so that $\abs{v}  \leq M$ for all $v \in \mathcal{H}$. From Lemma \ref{lem:tauuniform}, we know there is a function $r(t) \geq 0$ (depending on $\omega$, $M$ and $R$) such that $\lim_{t \to \infty} r(t) = 0$ and
\be
\sup_{\abs{x} \leq Rt} \;\; \sup_{\abs{v} \leq M}  \left| \frac{\tau(x,x + tv,\omega)}{t} - \bar q(v) \right| \leq r(t) \label{ttbound}
\ee
holds for all $t > 1$. This tells us that the travel time from any $x$ to $x + tv$ is approximately $t$. Specifically, for any $\abs{x} \leq Rt$ and $v \in \mathcal{H}$ we may choose a control $\alpha$ such that 
\[
X^\alpha_{x}(s) = x + tv
\]
for some time $s = \tau(x,x + tv,\omega)$ satisfying $\abs{s - t} \leq r(t) t$. If $s \leq t$, then the bound on $V$ implies
\[
\left| X^\alpha_{x}(t) - (x + tv) \right|  =  \left| X^\alpha_{x}(t) - X^\alpha_{x}(s)\right| \leq (1 + \norm{V}_\infty)\abs{t - s}.
\]
Similarly, if $s \geq t$, then we may extend the control $\alpha$ by setting $\alpha(r) = 0$ for $r \in [s,t]$. Using this modified control we obtain the same bound: $\left| X^\alpha_{x}(t) - (x + tv) \right|  \leq  (1 + \norm{V}_\infty)\abs{t - s} $. In either case, this shows that for all $\abs{x} \leq Rt$ and $v \in \mathcal{H}$, we may choose a control $\alpha \in \mathcal{A}_t$ such that 
\[
\left| X^\alpha_{x}(t) - (x + tv) \right| \leq (1 + \norm{V}_\infty) r(t)t
\]
Consequently, for all $\abs{x} \leq Rt$ and $z \in \partial \{x + W_t\} = x + t \mathcal{H}$ there must be a point $y \in \Gamma_{x,t}(\omega)$ such that $\abs{y - z} \leq (1 + \norm{V}_\infty) r(t)t$. Since $r(t) \to 0$ as $t \to \infty$, this proves (ii) for $z \in \partial \{x + W_t\}$. 

Finally, suppose that $z = x + tv$ is in the interior of $\{ x + W_t\} = \{ x + t W_1 \}$, with $\text{dist}(v,\mathcal{H}) \geq \epsilon_1$ and $\epsilon_1 \in (0,1)$. The set $W_t$ is convex. So, we may choose $s \in [0,1]$ and $v_1, v_2 \in \mathcal{H}$ such that $v = s v_1 + (1 - s) v_2$. 
Since $\text{dist}(v,\mathcal{H}) \geq \epsilon_1$ there must be $\epsilon_2 \in (\epsilon_1/(2M),1 - \epsilon_1/(2M))$ such that $s \in (\epsilon_2 , 1 - \epsilon_2)$. Since $v_1 \in \mathcal{H}$ there must be a point $y_1 \in \Gamma_{x,ts}(\omega)$ such that $\abs{y_1 - (x + t s v_1)} \leq \delta s t$ if $t \geq (\epsilon_2)^{-1} t_0(\delta,R,\omega)$. Since $v_2 \in \mathcal{H}$, there also must be a point $y_2 \in \Gamma_{y_1,t(1 -s)}(\omega)$ such that $\abs{y_2 - (y_1 + t (1 - s)v_2)} \leq \delta (1 -s)t$ if $t \geq (\epsilon_2)^{-1} t_0(\delta,R,\omega)$. Consequently, for all $\abs{x} \leq Rt$ and $t \geq (\epsilon_2)^{-1} t_0(\delta,R,\omega)$, there is $y_2 \in \Gamma_{x, t}(\omega)$ such that 
\br
\abs{y_2 - (x + t v)} & = & \abs{y_2 - (y_1 + t (1 - s)v_2) + y_1 - (x + t s v_1))}\leq \delta t
\er
holds. This completes the proof of (ii). \qed

\begin{corollary}\label{cor:asympdense} Let $\delta > 0$, $R > 0$, $h > 0$. For all $\omega \in \tilde \Omega$, there exists $\epsilon_0 = \epsilon_0(\delta,R,h,\omega) > 0$ such that if $\epsilon \in (0,\epsilon_0)$, the following hold:
\bi
\item[(i)] For all $\abs{x} \leq R$ and $t \geq h$, $\epsilon \Gamma_{x\epsilon^{-1},t\epsilon^{-1}}(\omega) \subset \{ x + W_{t(1 + \delta)}\}$.
\item[(ii)] For all $\abs{x} \leq R$, $t \geq h$, and all $z \in \{ x + W_t\}$ there is a point $y \in \epsilon \Gamma_{x\epsilon^{-1},t\epsilon^{-1}}(\omega)$ such that $\abs{y - z} \leq \delta t$.
\ei
\end{corollary}

{\bf Proof of Corollary \ref{cor:asympdense}:} Since $\epsilon W_{t \epsilon^{-1}} = W_t$, this is an immediate consequence of Lemma \ref{lem:asympdense}, with $\epsilon_0 = h^{-1} t_0(\delta, R,\omega)$. \qed

\vspace{0.2in} 

{\bf Proof of Theorem \ref{theo:homog}:} The solution $\bar u$ of (\ref{barueqn}) is given by the control representation
\[
\bar u(t,x) = \sup_{y \in \{ x + W_t\} } u_0(y),
\]
while for each $\epsilon > 0$,
\[
u^\epsilon(t,x) = \sup_{y \in \epsilon \Gamma_{x\epsilon^{-1},t\epsilon^{-1}}} u_0(y).
\]
Let $\delta > 0$, $h > 0$, and $R > 0$. Applying Corollary \ref{cor:asympdense} (i), we see that for $\abs{x} \leq R$, $t \geq h$, and $\epsilon \geq \epsilon_0(\delta,h,R,\omega)$,
\[
u^\epsilon(t,x) \leq \sup_{y \in \{ x + W_{t (1 + \delta)}\}} u_0(y) = \bar u(t(1 + \delta),x). 
\]
Therefore, because of the uniform continuity of $\bar u$ and the arbitrariness of $\delta > 0$, we conclude that
\be
\limsup_{\epsilon \to 0} \sup_{t \geq h} \sup_{\abs{x} \leq R} \left( u^{\epsilon}(t,x,\omega) - \bar u(t,x) \right) \leq 0. \label{limsupeps}
\ee

Now we use Corollary \ref{cor:asympdense} (ii) to obtain a lower bound, as follows. For all $\abs{x} \leq R$, $t \in [h,T]$, $z \in \{ x + W_t \}$, and $\epsilon \geq \epsilon_0(\delta,h,R,\omega)$, there is a point $y^* \in \epsilon \Gamma_{x\epsilon^{-1},t\epsilon^{-1}}(\omega)$ satisfying $\abs{y^* - z} \leq \delta T$. If $z \in \{ x + W_t \}$ is chosen so that
\[
\bar u(t,x) =  \sup_{y \in \{ x + W_t\} } u_0(y) = u_0(z) 
\]
then from (ii),
\br
u^\epsilon(t,x) = \sup_{y \in \epsilon \Gamma_{x\epsilon^{-1},t\epsilon^{-1}}} u_0(y)  \geq  u_0(y^*) 
& \geq & u_0(z) - \abs{u_0(z) - u_0(y^*)} \no \\
& \geq & \bar u(t,x) - \phi(\delta T) \no
\er
where $\phi$ is the modulus of continuity for $u_0(x)$. Therefore, since $\delta > 0$ was arbitrary, we conclude that
\be
\liminf_{\epsilon \to 0} \inf_{t \geq h} \inf_{\abs{x} \leq R} \left( u^{\epsilon}(t,x,\omega) - \bar u(t,x) \right) \geq 0. \label{liminfeps}
\ee
This proves that $u^\epsilon \to \bar u$ uniformly on compact sets in $(0,\infty) \times \R^d$. To obtain the locally uniform convergence down to time $t = 0$, we observe that for $x \in \R^d$, $u^\epsilon(t,x)$ satisfies
\br
\sup_{t \in [0,h]} \abs{u^\epsilon(t,x) - \bar u(t,x)} & \leq & \sup_{t \in [0,h]} \abs{u^\epsilon(t,x) - u_0(x)}  + \sup_{t \in [0,h]} \abs{\bar u(t,x) - u_0(x)} .
\er
Since $\abs{X^{\epsilon,\alpha}_x(t) - x} \leq t(1 + V_\infty)$, the first term on the right is bounded by
\be
\sup_{t \in [0,h]} \abs{u^\epsilon(t,x) - u_0(x)} \leq \sup_{\substack{y \in \R^d\\ \abs{y - x} \leq h (1 + V_\infty)}} \abs{u_0(y) - u_0(x)} \leq \phi(h (1 + V_\infty))
\ee
This and a similar bound on $\abs{\bar u(t,x,\omega) - u_0(x)}$ implies that
\be
\lim_{h \to 0} \left [ \limsup_{\epsilon \to 0} \sup_{\substack{x \in \R^d\\ t \in [0,h]}} \abs{u^\epsilon(t,x,\omega) - \bar u(t,x,\omega)}\right] = 0. \label{tlessh}
\ee

Finally, by combining (\ref{limsupeps}), (\ref{liminfeps}), and (\ref{tlessh}), we conclude that (\ref{homconv}) holds with probability one. The stated properties of $H$ follow immediately from the properties of $\bar q$ and the Corollary \ref{cor:Hp} at the end of the next section. \qed

\section{A Lower Bound on the Reachable Set}\label{sec:reachable}

In this section we prove an estimate on the growth of the reachable set $\Gamma_{x,t}(\omega)$ defined at \eqref{reachable}, under the assumption that $V$ is divergence-free. This estimate holds in all dimensions $d \geq 2$ and does not rely on conditions (\ref{gammaexp}) and (\ref{taubound}). Since this estimate makes no use of the statistical structure of the vector field $V$, we will suppress the dependence of $\Gamma_{x,t}$, $V(x)$, and $\tau(x,y)$ on $\omega \in \Omega$.

\begin{lemma}\label{lem:isoperim} Let $d \geq 2$. Assume that $V \in C^1(\R^d;\R^d)$ is divergence free. For $x \in \R^d$ and $t \geq 0$, let $\Gamma_{x,t} \subset \R^d$ be the set of points that are reachable at time $t$:
\[
\Gamma_{x,t}=\{ y \in \mathbb{R}^d\;|\;X^{\alpha}_{x}(t) = y, \quad \text{for some}\; \alpha \in {\mathcal{A}_t} \},
\]
and $X^{\alpha}_{x}(t)$ solves~\eqref{controlode} with $\eps=1$.
For all $x \in \R^d$ and $t \geq 0$
\be
\abs{\Gamma_{x,t}} \geq \omega_d t^{d}\label{isop}
\ee
where $\omega_d$ is the volume of the unit ball in dimension $d$.
\end{lemma}

{\bf Proof:} Consider the map $X(t,x,\kappa):[0,\infty) \times \R^d \times \overline{B_1(0)} \to \R^d$ defined by
\be
\frac{\partial}{\partial t} X(t,x,\kappa) =  -V(X(t,x,\kappa)) + \kappa , \quad t \geq 0 \label{kappacontr}
\ee
with $X(0,x,\kappa) = x$. For all $\kappa \in \overline{B_1(0)}$ and $x \in \R^d$, the matrix $M(t,x,\kappa) = D_\kappa X(t,x,\kappa)$ satisfies
\[
\frac{\partial}{\partial t} M_{i,j} =  - \sum_{\ell} \frac{\partial V^{i}}{\partial X_\ell}(X(t,x,\kappa)) M_{\ell,j} + \delta_{ij}, \quad 1 \leq i,j \leq d 
\]
with $M(0,x,\kappa) = 0$. Consequently, for $A(t) = -DV(X(t,x,\kappa))$, we have
\[
M(t,x,\kappa) = \int_0^t e^{\int_s^t A(r) \,dr} I \,ds = t I + \sum_{n=1}^\infty \int_0^t \frac{\left( \int_s^t A(r) \,dr \right)^n}{n!}\,ds.
\]
Therefore, by our assumptions on $V$, there is a constant $C_1$ such that $\norm{M(t,x,\kappa) - t I } \leq C_1 t^2$ for all $t \in (0,1)$, $x \in \R^d$ and $\kappa \in B_1(0)$, where $I:\R^d \to \R^d$ is the identity. Then using the contraction mapping theorem, one can show that for any small $\epsilon > 0$ there is a $t_0 = t_0(\epsilon) > 0$, small such for all $t \in (0,t_0]$ the image of the map $X(t,x,\cdot):B_{1}(0) \to \R^d$ contains a ball $B_{\delta }(X(t,x,0))$ with radius $\delta$ bounded below by $\delta   > (1 - \epsilon)t$. Since $\Gamma_{x,t}$ is obtained from the control problem with a larger set of controls (not just the constant controls), this ball must be contained in $\Gamma_{x,t}$, which shows that $\abs{\Gamma_{x,t}} \geq \omega_d((1 - \epsilon)t)^d$ for all $t \in (0,t_0]$ and $x \in \R^d$.

Now, let $t_1 \in (0,t_0)$ and define $t_k = k t_1$ for positive integers $k$. The analysis above shows that
\be
\Gamma_{x,t_{k+1}} \supset \bigcup_{y \in \Gamma_{x,t_{k}}} B_{\delta }(X(t_1,y,0)) \label{BM1}
\ee
must hold for $\delta = (1-\epsilon)t_1$. Since $V$ has zero divergence, the flow defined by (\ref{kappacontr}) with $\kappa = 0$ is volume-preserving. Hence, $\abs{\Gamma_{x,t}} \geq \abs{\Gamma_{x,s}}$ holds for all $t \geq s$ and
\be
\left| \bigcup_{y \in \Gamma_{x,t_{k}}} X(t_1,y,0)  \right| = \left| \Gamma_{x,t_{k}} \right| \label{BM2}
\ee
Now by applying the Brunn-Minkowski inequality to (\ref{BM1}) and using (\ref{BM2}), we obtain the bound
\[
\abs{\Gamma_{x,t_{k+1}}}^{1/d} \geq \delta \abs{\omega_d}^{1/d} + \left| \Gamma_{x,t_{k}} \right|^{1/d}. 
\]
Iterating this inequality yields
\[
\abs{\Gamma_{x,t_{k+1}}}^{1/d} \geq (k-1) \delta  \abs{\omega_d}^{1/d} + \left| \Gamma_{x,t_{1}} \right|^{1/d} \geq k \delta \abs{\omega_d}^{1/d}. 
\]
So, $\abs{\Gamma_{x,t_{k+1}}} \geq \omega_d (1 - \epsilon)^d (k t_1)^d$ holds for all $k \in \mathbb{N}$. Because $\abs{\Gamma_{x,t}} \geq \abs{\Gamma_{x,s}}$ holds for all $t \geq s$ and because $t_1$ and $\epsilon$ may be made arbitrarily small, this implies (\ref{isop}). \qed

\vspace{0.2in}

An immediate consequence of Lemma \ref{lem:isoperim} and Lemma \ref{lem:tauuniform} is the following simple estimate on the functions $\bar q$ and $\bar H$:

\begin{corollary} \label{cor:Hp} Suppose that $V$ is divergence-free and that conditions (\ref{gammaexp}) and (\ref{taubound}) hold. Then the function $\bar q$ defined by Lemma \ref{lem:subadd} satisfies
\be
\left| \left\{ z \;|\; \bar q(z) \leq 1 \right\} \right| \geq \omega_d
\ee
where $\omega_d$ is the volume of the unit ball in dimension $d \geq 2$. Moreover, if the law of $V(x,\omega)$ is equal to the law of $V(\Phi x,\omega)$ for any orthogonal transformation $\Phi$, then $\left\{ z \;|\; \bar q(z) \leq 1 \right\}  \supset B_1(0)$ and therefore $\bar H(p) \geq \abs{p}$ for all $p \in \R^d$.
\end{corollary}

\section{An Upper Bound on the Travel Times  in Dimension $d = 2$} \label{sec:ttbound}

In this section we prove Theorem \ref{thm:twodim} by showing that the conditions (\ref{gammaexp}) and (\ref{taubound}) hold for a large class of vector fields in dimension $d = 2$. If $V \in C^1(\R^2;\R^2)$ is divergence-free, then there is a stream function $\F:\R^2 \to \R$ such that
\[
V(x) = \nabla^\perp \F(x) = (- \partial_2 \F, \partial_1 \F).
\]
Some of the analysis of this section does not use the statistical structure of $V$ and $\F$ so we suppress the dependence of $V$, $\F$, and $\tau$ on $\omega \in \Omega$. Our main estimate is the following lemma which relates the travel times $\tau(x,y)$ to the growth of the stream function:

\begin{lemma} \label{prop:travtimebound} There are constants $C_1$ and $C_2$ such that if $M > 0$, $K > 1$, $z \in \R^d$, and 
\be
\abs{\F(x) - \F(z)} \leq K \label{HKassump}
\ee
holds for all $x$ satisfying $\abs{x - z} \leq 3M + 5 K$, then $\tau(x,y) \leq C_1 K+C_2 \abs{x - y}$ holds for all $x,y \in \overline{B_M(z)}$. 
\end{lemma}

To prove this bound, we introduce a modified control problem associated with a modified stream function $\hat \F$.

\begin{lemma}\label{Hp} Suppose $\F$, $M > 0$, $K > 1$, and $z \in \R^d$ are such that \eqref{HKassump} holds for all $\abs{x - z} \leq 3M + 5 K$. Then we can modify  $\F$ by adding a non-negative function $\phi$ so that the modified function
\[
\hat \F=\F+\phi
\]
satisfies:
\begin{itemize}
\item  $\hat \F(x) = \F(x)$ if  $\abs{x - z} \leq M$, 
\item $K +1/3 \leq \hat \F(x) - \hat \F(z) \leq 3 K+1/3+M$, for $ M+1+ 4K \leq \abs{x - z} \leq 3M+1+4K$, 

\item $|\nabla \phi|\leq1/2$. 
\end{itemize}
\end{lemma}
{\bf Proof of Lemma \ref{Hp}:} 
Let $\rho(s) \in C^\infty(\mathbb{R^+})$ be a non-decreasing function that satisfies
\[
\rho(s)=\begin{cases}
0,& s \leq M,\\
s/2 -M/2-1/6,& s \geq M+1,\\
\end{cases}
\]
and $0 \leq \rho(s) \leq 1/3$ when $M \leq s \leq M+1$.
We can choose $\rho$ so that $0 \leq d\rho(s)/ds \leq 1/2$.
Then we define $\phi(x) =   \rho(\abs{x - z})$.
\qed

The modified stream function $\hat \F$ depends on $z$. For $z$ fixed and $\hat \F$ defined in this way, let $\hat V(x)=\nabla^{\perp} \hat \F$ be the vector field associated with the modified stream function, and define an auxiliary control problem
\begin{equation}\label{c2}
\frac{d}{dt} Y^{\hat \alpha}_x(t)= -\hat V(Y^{\hat \alpha}_x(t))+ \hat \alpha(t),\quad Y^{\hat \alpha}_x(0) = x \in\mathbb{R}^2,
\end{equation}
where the control $\hat \alpha$ satisfies the constraint $|\hat \alpha|\leq 1/2$:
\[
\hat \alpha \in \hat {\mathcal{A}_t} = \{ \hat \alpha \in L^\infty([0,t];\R^d) \;:\; \| \hat \alpha \|_\infty \leq 1/2 \} \subset \mathcal{A}_t.
\]
Travel times for the unmodified control problem are bounded by travel times for the modified problem:

\begin{lemma}\label{taup} For any solution $Y^{\hat \alpha}_x(t)$ with control $\hat \alpha \in \hat {\mathcal{A}}_t$ of the system~\eqref{c2}, there is a solution $X^\alpha_x(t)$ with control $\alpha \in \hat {\mathcal{A}}_t$ of the 
system~\eqref{controlode} such that $X^\alpha_x(s) =Y^{\hat \alpha}_{x}(s)$ for all $s \in [0,t]$. Thus, 
\be
\tau(x,y) \leq  \hat \tau(x,y), \label{tauhattau}
\ee
where 
\[
\hat \tau(x,y)= \inf\{t \geq 0 \;|\;Y_x^{\hat \alpha}(t) = y, \quad \text{for some} \;\; \hat \alpha \in \hat {\mathcal{A}_t} \}.
\]

\end{lemma}
{\bf Proof of Lemma \ref{taup}:}  Given a control $\hat \alpha \in \hat {\mathcal{A}}_t$, we set $\alpha(s)=\nabla^\perp \phi(Y^{\hat \alpha}_x(s)) +\hat \alpha(s)$ for $s \in [0,t]$ and solve (\ref{controlode}) with this control. By definition of $\phi$, $\alpha \in \mathcal{A}_t$. Also, $X^{\alpha}_x(s) = Y^{\hat \alpha}_x(s)$ for all $s \in [0,t]$. If $\hat \tau(x,y)$ is infinite, the bound (\ref{tauhattau}) holds trivially.  If $\hat \tau(x,y) < \infty$, there is a time $t$ and a control $\hat \alpha \in \hat {\mathcal{A}_t}$ such that $Y^{\hat \alpha}_x(t) = y$, so by choosing $\alpha$ as we have just described, we obtain (\ref{tauhattau}). \qed

\vspace{0.2in}

{\bf Proof of Lemma \ref{prop:travtimebound}:} Suppose that $M > 0$, $K > 1$, and $z \in \R^d$ satisfy the hypotheses of the Lemma. Let us choose $\phi(x) = \rho(\abs{x - z})$ as in Lemma~\ref{Hp}, and $\hat \F = \F + \phi$.  First we will show that if $x_0, y_0 \in \overline{B_M(z)}$, then 
\begin{equation}\label{dTaup}
\tau(x_0,y_0) \leq C_1 K+C_2 M
\end{equation} 
must hold for some constants $C_1$ and $C_2$ that are independent of $x_0, y_0, K, M, z$. Lemma \ref{Hp} implies that for all $h$ bounded by $K +1/3 \leq h \leq 3 K+1/3+M$
there exists a connected component 
$S^0_h$ of the set 
$S_h =\{ x \in \R^d \;|\;  \hat \F(x) - \hat \F(z) \leq h\}$ satisfying $\overline{B_{M}(z)} \subset S^0_h \subset B_{R}(z)$,  
where $R=3M+1+4K \leq 3M + 5K$. 

We will obtain~\eqref{dTaup} from a similar estimate on $\hat \tau(x_0,y_0)$. 
To estimate $\hat \tau(x_0,y_0)$ we will estimate the time required to move from $x_0$ to any point in the set $\partial S^0_h$, then the time to reach 
$y_0$ from some point in $\partial S^0_h$.

We first estimate the time to reach the boundary of $B_{R}(z)$, starting from $x_0$. For each $t\geq 0$ we define the (modified) reachable set
\[
\hat \Gamma_{x_0,t}=\{ y \in \mathbb{R}^2\;|\;Y^{\hat \alpha}_{x_0}(t) = y, \quad \text{for some}\;\hat \alpha \in \hat {\mathcal{A}_t} \},
\]
and $Y^{\hat \alpha}_{x}(t)$ solves~\eqref{c2}.
From Lemma \ref{lem:isoperim} it follows that in dimension $d = 2$,
\be
\abs{\hat \Gamma_{x_0,t}}  \geq   \frac{\pi t^2}{4}, \quad \forall \;\; t \geq 0. \label{gammalower}
\ee
holds for all $t > 0$. The factor $1/4$ here comes from the fact that the controls in $\hat{\mathcal{A}_t}$ defining $Y^{\hat \alpha}$ must satisfy $\abs{\hat \alpha} \leq 1/2$, rather than $\abs{\alpha} \leq 1$. If $\hat \Gamma_{t}^+ = \cup_{s \in [0,t]} \hat \Gamma_{x_o,s}$, the lower bound  (\ref{gammalower}) implies that $\hat \Gamma_t^+ \cap \partial B_{R}(z) \neq \emptyset$ for times  $t > 2 R$. Thus the curve $\gamma=\partial S^0_h$ also intersects $\hat \Gamma_t^+$ at times $t > 2R$, and for all $h \in [K +1/3,  3 K+1/3+M]$.

The curve $\gamma$ is an integral curve of a solution of~\eqref{c2} with $\hat \alpha \equiv 0$, thus $\gamma \subset \partial S_h$ for some $h$. Let us estimate the time 
\[
T=\int_{\gamma} \frac{1}{|\hat{V}|}  ds 
\]
required to traverse this curve. Since $\F \in C^2$ and $\phi \in C^{\infty}$ we have $\hat{\F} \in C^2$. Applying the Sard's lemma we conclude that the set $\mathcal R$ of regular values of $\hat{\F}$ has full measure. 
  The implicit function theorem guarantees $\partial S_h$ is a finite union of $C^2$ compact manifolds without boundary for each $h\in \mathcal R$. $\mathcal R$ is also open because $\hat{\F}\in C^2(S)$.
For any $\delta>0$ by the co-area formula 
\[
\int_{h_1}^{h_2} \left( \int_{\partial S_h} \frac{1}{|\nabla \hat{\F}|+\delta}  ds \right) dh =  \int_{S_{h_2} \setminus S_{h_1}} \frac{|\nabla \hat{\F}|}{|\nabla \hat{\F}|+\delta} 
dx dy \leq |S_{h_2}| -   |S_{h_1}|,
\]
where $h_1 \leq h_2$, and $|S_h|$ denotes the area of $S_h$. Noting that $| \hat V|=|\nabla \hat{\F}|$ and using the monotone convergence theorem we obtain
\[
\int_{h_1}^{h_2} \left( \int_{\partial S_h} \frac{1}{|\hat{V}|}  ds \right) dh =  |S_{h_2}| -   |S_{h_1}|.
\]
Since $|S_{K +1/3}| \geq \pi M^2$ and $|S_{3 K+1/3+M}| \leq  \pi R^2$ the last identity implies
\[
\int_{K+1/3}^{3 K +1/3 +M} \left( \int_{\partial S_h} \frac{1}{|\hat{V}|}  ds \right) dh \leq \pi \left( R^2-M^2 \right).
\]
Thus there exists $h \in \mathcal{R} \cap [K +1/3,  3 K+1/3+M]$ so that
\begin{equation}\label{chooseh}
T = \int_{\partial S^0_h} \frac{1}{|\hat{V}|}  ds \leq  \int_{\partial S_h} \frac{1}{|\hat{V}|}  ds \leq\pi \left( R^2-M^2 \right)/\left(2K+M \right)\leq  4 \pi R \leq 13 R.
\end{equation}
Thus at all times $t> (13 + 2) R $, $\hat \Gamma^+_t$ contains the $C^2$ curve $\gamma=\partial S^0_h \subset S_h$, where $h$ is chosen, so that~\eqref{chooseh} 
holds.

We complete the proof of (\ref{dTaup}) by estimating the time needed for $\hat \Gamma^+_t$ to contain the entire
set $S^0_h$. If we replace $\hat V$ with $-\hat V$, then starting from any point $y_1 \in S^0_h$ at time $t=0$, we may apply Lemma \ref{lem:isoperim}, as above, to show that there is a control $\beta \in \hat {\mathcal{A}_t}$ for which the controlled path reaches $\partial S^0_h$ at some time $t^* \leq 2R$. Therefore by reversing time with this control, we conclude that for any point $y_1 \in \overline{B_M(z)} \subset S^0_h$ there is a control $\alpha(s) = -\beta(t^* - s) \in \hat{\mathcal{A}}_{t_*}$ and a point $y_2 \in \partial S^0_h$ for which $Y_{y_2}^\alpha(t^*) = y_1$. 
Thus $y_1 \in \hat \Gamma^+_t$ for all $t > (13 + 2 +2) R = 17R$. Using our definition of $R$ we 
obtain the bound $\hat \tau(x_0,y_0)\leq C_1 K + C_2 M$ for all $x_0, y_0 \in \overline{B_M(z)}$. Now \eqref{dTaup} follows immediately from this and Lemma \ref{taup}. 

Now, using (\ref{dTaup}), we show that
\be
\tau(x_0,y_0) \leq C_1' K + C_2' \abs{x_0 - y_0} \label{taux0y0new}
\ee
holds for all $x_0, y_0 \in \overline{B_M(z)}$, with constants $C_1'$ and $C_2'$ that are independent of $M,K,z,x_0, y_0$. If $\abs{x_0 - y_0} \geq M/4$, then it is immediate from (\ref{dTaup}) that $\tau(x_0,y_0) \leq C_1 K + C_2 M \leq C_1 K + 4 C_2 \abs{x_0 - y_0}$. In this case the estimate (\ref{taux0y0new}) holds with $C_1' = C_1$ and $C_2' = 4 C_2$.  So, let us suppose that $x_0, y_0 \in \overline{B_M(z)}$, but $\abs{x_0 - y_0} \leq M/4$. If $M < 5K$, then $\tau(x_0,y_0) \leq C_1 K + C_2 M \leq (C_1 + 5C_2) K$. So, (\ref{taux0y0new}) certainly holds with $C_1' = (C_1 + 5C_2)$ and $C_2' = 0$. Therefore, it suffices to assume that $\abs{x_0 - y_0} \leq M/4$ and $M \geq 5K$. In this case, let us define $z' = (x_0 + y_0)/2$, $M' = \abs{x_0 - y_0} \leq M/4$, and $K' = 2K$. Then it is easy to see that
\[
B_{3M' + 5K'}(z') \subset B_{3M + 5K}(z)
\]
since $3M' + 5 K' < 2M + 5 K$ and $z' \in B_{M}(z)$. Therefore, $|\F(x) - \F(z')| \leq K'$ holds for all $\abs{x - z'} \leq 3M' + 5K'$. Thus, we may apply the bound (\ref{dTaup}) with $M'$, $K'$ and $z'$, to conclude that $\tau(x_0,y_0) \leq  C_1 K'+C_2 M' = 2 C_1 K + C_2 \abs{x_0 - y_0}$. In this case (\ref{taux0y0new}) holds with $C_1' = 2 C_1$ and $C_2' = C_2$. Finally, by combining each of these cases, we see that (\ref{taux0y0new}) holds for all $x_0, y_0 \in \overline{B_M(z)}$ with $C_1' = 2C_1 + 5 C_2$ and $C_2' = 4C_2$. \qed

\vspace{0.2in}

In proving Theorem \ref{thm:twodim}, we will make use of the following lemma which shows that the stream function grows no more than sublinearly, if $\expE[V] = 0$:

\begin{lemma} \label{lem:streamgrowth} Let $V:\R^2 \times \Omega \to \R^2$ be a stationary random vector field, uniformly bounded and divergence-free. Let $\F(x,\omega): \R^2 \times \Omega \to \R$ be a stream function: $V = \nabla^\perp \F = (- \partial_{x_2} \F, \partial_{x_1} \F)$. If $\expE[V_i(0,\omega)] = 0$, for $i=1,2$, then $\F$ must satisfy 
\be
\lim_{r \to \infty} \frac{1}{r} \sup_{\abs{x} \leq r} \abs{\F(x,\omega)} = 0 \label{streamsub}
\ee
with probability one.
\end{lemma}
{\bf Proof of Lemma \ref{lem:streamgrowth}:} Since $V$ is uniformly bounded, $\F(x,\omega)$ must be uniformly Lipschitz continuous. Therefore, it suffices to prove that with probability one,
\be
\lim_{r \to \infty} \frac{1}{r} \abs{\F(r z,\omega)} = 0 \label{streamconv}
\ee
holds for any fixed vector $z \in \R^2$ with rational coordinates (i.e. $z \in \mathbb{Q}^2$).

Without loss of generality, we may assume that $\F(0,\omega) = 0$ almost surely.  So, for any curve $p:[0,r] \to \R^2$ we have
\[
\F(p(r),\omega)  =  \int_0^r p' \cdot \nabla \F \,dt =\int_0^r p_1'(t) V_2(p(t),\omega) - p_2'(t)V_1(p(t),\omega) \,dt.
\]
For a given $z \in \mathbb{Q}^2$, let $p(s) = z s$. Therefore,
\[
\lim_{r \to \infty} \frac{1}{r} \F(r z,\omega) = \lim_{r \to \infty} \frac{1}{r} \int_0^r \theta(t,\omega) \,dt
\]
where $\theta(t,\omega) = z_1 V_2(z t,\omega) - z_2 V_1(z t,\omega)$. Observe that $\theta(t + h,\omega) = \theta(t,\pi_{zh} \omega)$, so that $\theta(t,\omega)$ is statistically stationary in $t$.  Therefore, the ergodic theorem implies that there is a random variable $\eta \in L^1(\Omega)$ such that the limit
\be
\lim_{r \to \infty} \frac{1}{r} \int_0^r \theta(t,\omega) \,dt = \eta(\omega)
\ee
holds with probability one. We claim that $\eta \equiv 0$. This follows from the fact that $\eta$ must be invariant under the action of $\pi_{y}$ for all $y \in \R^2$. So see this, observe that the stationarity of $V$ implies that $\int_0^r \theta(t,\pi_y \omega) \,dt  =  \F(zr + y,\omega) - \F(y,\omega)$. Because $\abs{(\F(zr + y,\omega) - \F(y,\omega)) - \F(zr,\omega)} \leq 2 K y$ we then have
\br
\eta(\pi_y \omega) & = & \lim_{r \to \infty} \frac{1}{r}  \left(\F(zr + y,\omega) - \F(y,\omega) \right)  \no \\
& = & \lim_{r \to \infty} \frac{1}{r}  \F(zr,\omega) + \lim_{r \to \infty} \frac{1}{r}  \left((\F(zr + y,\omega) - \F(y,\omega)) - \F(zr,\omega) \right) \no \\
& = & \eta(\omega). \no
\er
Thus, $\eta$ is invariant under $\pi_y$ for all $y \in \R^2$, so $\eta$ must be a constant, and therefore, $\eta = \expE[\eta] = \expE[ z^\perp \cdot V(0,\omega)] = 0$. This establishes (\ref{streamconv}). \qed

\vspace{0.2in}

{\bf Proof of Theorem \ref{thm:twodim}:} Our proof consists of verifying conditions (\ref{gammaexp}) and (\ref{taubound}) in Theorem \ref{theo:homog}.  First, let us verify condition (\ref{taubound}), using the assumption that $\expE[V] = 0$. By Lemma \ref{prop:travtimebound}, with $K = M_0 = M$ and $z = 0$, we know there is a constant $C_3$ such that if $\abs{\F(x,\omega) - \F(0,\omega)} \leq M_0$ holds for all $\abs{x} \leq 8M_0$, for some $M_0 > 1$, then $\gamma(M_0,\omega) \leq C_3 M_0$. From Lemma \ref{lem:streamgrowth}, we know that with probability one there is a random variable $M_0(\omega) < \infty$ such that if $M \geq M_0(\omega)$ then $\abs{\F(x,\omega) - \F(0,\omega)} \leq M$ holds for all $\abs{x} \leq 8M$. Consequently,
\be
\gamma(M,\omega) \leq C_3 M, \quad \forall \;M \geq M_0(\omega)
\ee
holds with probability one.  This implies condition (\ref{taubound}).

Now we verify condition condition (\ref{gammaexp}), using assumption (\ref{streamgrowthcond}). Fix $M \geq 2$. For any $K > 3M/5$, Lemma \ref{prop:travtimebound} implies that
\br
\Pm \left( \sup_{x, y \in \overline{B_M(0)} } \tau(x,y) > C_1 K + C_2 M \right) &\leq & \Pm \left( \sup_{\abs{x} \leq 3 M + 5K} \abs{\F(x) - \F(0)} \geq K \right)  \no \\
& = & \Pm \left( \sup_{\abs{x} \leq R} \abs{\F(x) - \F(0)} \geq \frac{R}{5} - \frac{3M}{5} \right)  
\er
where $R = 3M + 5K$. If $R > 900M/5$, then $R/5 - 3M/5 > R/6$, so in this case we have
\be
\Pm \left( \sup_{x, y \in \overline{B_M(0)} } \tau(x,y) > C_1 K + C_2 M \right) \leq \Pm \left( \sup_{\abs{x} \leq R} \abs{\F(x) - \F(0)} \geq \frac{R}{6} \right).  
\ee
By assumption (\ref{streamgrowthcond}), this last quantity is integrable in $R$, so
\br
\expE \left[  \sup_{x, y \in \overline{B_M(0)} } \tau(x,y) \right]& = & \int_0^\infty \Pm \left( \sup_{x, y \in \overline{B_M(0)} } \tau(x,y) > r \right) \,dr \no \\
& \leq & C_4 M + C_5 \int_{0}^\infty \Pm \left( \sup_{\abs{x} \leq R} \abs{\F(x) - \F(0)} \geq \frac{R}{6} \right)  \,dR \no \\
& \leq & C_4 M + C_6.
\er
This establishes condition (\ref{gammaexp}). \qed

\vspace{0.2in}

{\bf Proof of Proposition \ref{prop:Fstat}:} For any $R > 1$, we may cover the set $B_R(0)$ by $O(R^2)$ balls of radius $1$. That is, there is a constant $C$, independent of $R$, such that for each $R > 1$ there is a set of points $\{ x_i \}_{i=1}^N$ with $N \leq C R^2$ such that $B_R(0) \subset \bigcup_{i=1}^N B_1(x_i)$. Therefore, since $\F$ is stationary,
\br
\Pm \left( \sup_{\abs{x} \leq R} \abs{\F(x) - \F(0)} \geq R/6 \right) & \leq & \Pm \left( \sup_{\abs{x} \leq R} \abs{\F(x) } \geq R/12 \right)  \no \\
 & \leq & \sum_{i=1}^N \Pm \left( \sup_{x \in B_1(x_i)} \abs{\F(x)} \geq R/12	 \right) \no \\
& \leq & C R^2 \Pm \left( \sup_{\abs{x} \leq 1} \abs{\F(x)} \geq R/12 \right) \label{coverbound}
\er
holds for all $R > 1$. Consequently
\br
\int_{1}^\infty \Pm \left( \sup_{\abs{x} \leq R} \abs{\F(x) - \F(0)} \geq \frac{R}{6} \right)  \,dR  & \leq &\int_1^\infty C R^2 \Pm \left( \sup_{\abs{x} \leq 1} \abs{\F(x)} \geq R/12 \right) \,dR \no \\
& \leq & C' \expE \left[ \left| \sup_{\abs{x} \leq 1} \abs{\F(x)} \right|^{3}  \right] 
\er
Since $V$ is uniformly bounded by $V_\infty$, then $\abs{\nabla \F} \leq V_\infty$. Therefore, the last expectation is bounded by
\[
\expE \left[ \left| \sup_{\abs{x} \leq 1} \abs{\F(x)} \right|^{3}  \right] \leq \expE \left[ \left(\abs{\F(0)} + V_\infty  \right)^{3}  \right], 
\]
which is finite, by assumption (\ref{Fmomentbound}). Thus condition (\ref{streamgrowthcond}) holds. \qed

\end{document}